\documentclass[final]{elsart}
\usepackage{amsmath,amssymb,latexsym}
\usepackage[final]{graphicx}

\newcommand{\supp}{\mathop{\rm supp}}
\newcommand{\field}[1]{\mathbb{#1}}
\newcommand{\R}{\field{R}}

\newcommand{\N}{\field{N}}
\newcommand{\C}{\field{C}}
\newcommand{\T}{\field{T}}
\renewcommand{\P}{\field{P}}
\newcommand{\F}{{\mathcal F}}

\newcommand{\const}{{\rm const}}

\renewcommand{\Im}{\mathop{\rm Im}}

\newcommand{\vt}[1]{\boldsymbol{#1}}

\newenvironment{proof}%
{\rm \trivlist \item[\hskip \labelsep{\bf Proof. }]}%
{\hspace*{\fill}$\Box$\endtrivlist}
{\rm \trivlist \item[\hskip \labelsep{\bf Proof}]}%
{\hspace*{\fill}$\Box$\endtrivlist}
\numberwithin{equation}{section}

\begin{document}

\begin{frontmatter}
%%%%%%%%%%%%%%%%%%% Frontmatter %%%%%%%%%%%%%%%%%%%%%%

\title{Electrostatic models for zeros of polynomials: old, new, and some open problems}

\author[Madrid]{F.\ Marcell\'{a}n\corauthref{cor}}
\author[Almeria,Granada]{, A. Mart\'{\i}nez-Finkelshtein}
\author[Almeria]{and P. Mart\'{\i}nez-Gonz\'{a}lez}
\address[Madrid]{Departamento de Matem\'{a}ticas, Escuela Polit\'{e}cnica Superior,
Universidad Carlos III, Avenida de la Universidad, 30, 28911
Legan\'{e}s-Madrid, Spain}
\address[Almeria]{Departamento de Estad\'{\i}stica y Matem\'{a}tica Aplicada,
Universidad de Almer\'{\i}a, \\ La Ca\~{n}ada, 04120 Almer\'{\i}a,
Spain}
\address[Granada]{Instituto Carlos I de F\'{\i}sica Te\'{o}rica y
Computacional, Universidad de Granada, 18071 Granada, Spain}
\corauth[cor]{Corresponding author. E-mail:
\texttt{pacomarc@ing.uc3m.es}}

{\sc Dedicated to Nico Temme on the occasion of his 65th birthday.}

\begin{abstract}
We give a  survey  concerning  both very classical and recent
results on the electrostatic interpretation of the zeros of some
well-known families of polynomials, and the interplay between these
models and the asymptotic distribution of their zeros when the
degree of the polynomials tends to infinity. The leading role is
played by the differential equation satisfied by these polynomials.
Some new developments, applications and open problems are presented.
\end{abstract}


%\amsclassification{26C10, 30E15, 31A35, 33C45, 34C10, 42C05}

% {\bf AMS (MOS) Subject Classification:}

\end{frontmatter}
% ------------------------------------------------------------

\section{Introduction}

The electrostatic interpretation of the zeros of the classical
orthogonal polynomials is probably one of the most elegant results
in the theory of special functions, linked in the first instance to
such a distinguished name as Stieltjes (although studied also by
B\^{o}cher, Heine, Van Vleck, and Polya). Although this topic has
remained ``dormant'' for almost a century, a renewed interest has
appeared recently, partially motivated by the connection of this
topic with modern powerful techniques from the theory of logarithmic
potentials,  as well as by an increasing interest on the study of
new classes of special functions with roots in Physics,
Combinatorics, Number Theory, etc.

This paper is a short and light survey on these topics. Our
intention is to provide not a comprehensive list of results, but
more to convey the ``flavor'' of techniques and ideas behind the
electrostatic model and its direct implications in the asymptotic
theory. With this purpose, in Section \ref{sec2} we derive the
classical electrostatic interpretation of the zeros of Jacobi
polynomials, found originally by Stieltjes, and we describe some of
its new and still-in-progress direct generalizations (Section
\ref{sec3}). Next we discuss a recent progress by Ismail in
connecting electrostatics with orthogonality in Section
\ref{secOrthogonality}. Seeking for a model for complex zeros in
Section \ref{secMaxMin} we formulate an alternative model, not
constrained to the real line, as well as discuss some further
generalizations and applications of the electrostatic model in
Section \ref{secFurther}. We conclude the paper with some examples
of applications of these models to the study of the asymptotic
distribution of zeros in the semiclassical limit.

\section{Electrostatic model for classical orthogonal polynomials}
\label{sec2}

Let us begin with a classical and  very well-known result, found by
Stieltjes in 1885,  about  the electrostatic interpretation of the
zeros of Jacobi polynomials $P_n^{(\alpha,\beta)}$. The definitions
and basic properties of these polynomials can be found in Chapter IV
of the classical monograph  by  Szeg\H{o} \cite{Szego:75}; we will
mention here only those that will be needed further.

Jacobi polynomials can be given explicitly by
$$
P_n^{(\alpha,\beta)} (x)=2^{-n} \sum_{k=0}^n \binom{n+\alpha}{n-k}
\binom{n+\beta}{k}(x-1)^k (x+1)^{n-k}\,.
$$
Nevertheless, this formula is not the most useful one, at least for
what we are looking for. Another equivalent characterization is the
Rodrigues formula
\begin{equation}\label{rodrigues}
P_n^{(\alpha,\beta)} (x)=\frac{1}{2^n n!} (x-1)^{-\alpha}
(x+1)^{-\beta} \left( \frac{d}{dx} \right)^n \left[ (x-1)^{n+\alpha}
(x+1)^{n+\beta}\right]
\end{equation}
(cf.\ Section 4.3 of \cite{Szego:75}).  In particular, these
expressions show that Jacobi polynomials $P_n^{(\alpha,\beta)}$ are
analytic functions of the parameters $\alpha, \beta \in \C$ and that
$\deg P_n^{(\alpha,\beta)} \leq n$.

A third characterization of these polynomials, which will play the
leading role in the sequel, is that they are the only polynomial
solutions (up to a  constant factor) of the linear differential
equation
\begin{equation}\label{EDJac2}
y''(x)+\left(\frac{\alpha+1}{x-1}+\frac{\beta+1}{x+1}\right)y'(x)-
\frac{\lambda_{n}}{(x^2-1)}\, y(x)=0\,,
\end{equation}
where $\lambda_{n}=n(n+\alpha+\beta+1)$. This is a second order
linear differential equation of hypergeometric type.

The setting is considered ``classical'' when $ \alpha >-1$ and
$\beta >-1$. In this case, as we see from (\ref{EDJac2}), the
residues of the rational coefficient of $y'$ in the differential
equation are both positive. Moreover, we know that with this
assumption the Jacobi polynomials are orthogonal on $[-1,1]$ with
respect the weight function $(1-x)^{\alpha} (1+x)^{\beta}$, i. e.
\begin{equation}\label{orthogonalityForJacobi}
\int_{-1}^1 P_n^{(\alpha,\beta)}(x)x^k\, (1-x)^{\alpha}
(1+x)^{\beta}\, dx=0\,, \,  \quad k=0, 1, \dots, n-1\,.
\end{equation}
If both $ \alpha >-1$ and $\beta >-1$, then the weight is integrable
on $[-1,1]$, and a well known consequence of
\eqref{orthogonalityForJacobi} is that all the zeros of
$P_n^{(\alpha,\beta)}$ are simple and lie on the interval $(-1,1)$.
Apparently, it was Stieltjes \cite{Stieltjes:885} who  observed
first that it is possible to give a nice interpretation of the
location of these zeros as follows. Put two positive fixed charges
of mass $(\beta+1)/2$ and $(\alpha+1)/2$ at $-1$ and $+1$,
respectively, and allow $n$ positive \emph{unit} charges $X=\{x_1,
\dots, x_m \}$ to move freely in $(-1,1)$. If the interaction obeys
the logarithmic potential law (that is, the force is inversely
proportional to the relative distance), then in order to find the
total energy $E(X)$ of this system we have to add to the energy of
the mutual interaction of these charges,
$$
E_{\rm mutual}(X)=-\sum_{1\leq k < j \leq n}\ln|x_k-x_j|\,,
$$
the component given by the ``external field'' $\varphi(x)$ created
by the fixed charges,
\begin{equation}\label{external_field_Jacobi}
\varphi(x)=-\frac{\beta+1}{2}\, \ln|x+1|- \frac{\alpha+1}{2}\,
\ln|x-1|\,.
\end{equation}
In other words, the total energy is
\begin{equation}\label{energia}
E(X)=E_{\rm mutual}(X)+\sum_{k=1}^n\varphi(x_k) \,.
\end{equation}
There is a unique configuration $X^*=\{x_1^*, \dots, x_m^* \}$, $-1<
x_1^*<  x_2^* < \dots < x_n^* <1$, providing the (strict) global
minimum of $E(X)$ in $[-1,1]^n$, corresponding to the unique
equilibrium position for our free charges.

There uniqueness of the global minimum is not obvious, but there is
an elegant proof in \cite{Szego:75}, based on the inequality between
the arithmetic and geometric means. However, it is also a
consequence of the following statement (Stieltjes' theorem):
\emph{points $x_k^*$ are precisely the zeros of the polynomial
$P_n^{(\alpha,\beta)}$.} The proof is quite straightforward and uses
the differential equation \eqref{EDJac2}. Indeed, since $E(X) \to
+\infty$ as $X $ approaches any boundary point of $[-1,1]^n$, we
conclude that $X^* \subset (-1,1)^n$. Then $X^*$ is also a critical
point for the energy functional $E(X)$, and as a consequence, every
$x_k^*$ is a critical point of $E$ as a function of $x_k$ (fixing
the rest of $x_j^*$'s). Thus, the following necessary conditions
must be satisfied:
\begin{equation}\label{equations_critical}
\frac{\partial}{\partial
x_k}\,E(X)\big|_{X=X^*}=\frac{\partial}{\partial x_k}\,
E_{mutual}(X)\big|_{X=X^*}+\varphi'(x_k)=0\quad \text{for} \quad
k=1,2,\dots,n\,.
\end{equation}
Let us consider the monic polynomial vanishing at $x_k^*$'s:
$y(x)=\prod_{j=1}^{n}(x-x_j^*)$; it is easy to check that
$$
\frac{\partial}{\partial x_k}\, E_{mutual}(X)\big|_{X=X^*}=-
\sum_{1\leq j \leq n,\,j\neq
k}\frac{1}{x_k^*-x_j^*}=-\frac{y''(x_k^*)}{2 \, y'(x_k^*)}\,,
$$
so that \eqref{equations_critical} implies
\begin{equation}\label{equations_critical2}
y''(x) - 2 \varphi'(x) y'(x)=0\quad \text{for }  x \in X^*\,.
\end{equation}
For the external field \eqref{external_field_Jacobi} it means that
\begin{equation}\label{ceros}
y''(x)+\left(\frac{\beta+1}{x+1}+
\frac{\alpha+1}{x-1}\right)y'(x)=0\quad \text{for }   x \in X^*\,.
\end{equation}
However, the left hand side in \eqref{ceros} is a rational function
of the type $[n/2]$ with poles at $\pm 1$, so, up to a constant
factor, it is equal to $y(x)/(x^2-1)$, which yields (\ref{EDJac2})
and shows that $y(x)=\const \, P_n^{(\alpha,\beta)}(x)$.

Observe that we started from the strongest requirement of the
\emph{global} minimum of $E(X)$ (that we refer to as a \emph{stable
equilibrium}) which implied that $X^*$ is also a point of
\emph{local} minimum, and in consequence, we have a
\emph{``Nash-type'' equilibrium} (each function $E(x_k)=E(x_1^*,
\dots, x_{k-1}^*, x_k, x_{k+1}^*, \dots, x_n^*)$ attains its minimum
at $x_k=x_k^*$), yielding that each $x_k^*$ is a \emph{critical
point} of $E(x_1^*, \dots, x_{k-1}^*, x_k, x_{k+1}^*, \dots,
x_n^*)$. Precisely this last statement turned out to be equivalent
to the characterizing differential equation (\ref{EDJac2}). The
reader should take note of this clear hierarchy of equilibria in
this problem.

Interesting enough, in the situation we are considering all the
equilibria are equivalent. In order to prove it we can consider the
Hessian matrix
\begin{equation}\label{hesssian}
H=(h_{ij})\,,  \quad h_{ij}=\frac{\partial^2 E(X)}{\partial x_i
\partial x_j}
\end{equation}
of $E(X)$. Straightforward computations show that
$$ h_{ij} =
  \begin{cases}
    -2(x_i - x_j)^{-2} & \text{if } i \neq j, \\
   \displaystyle{ \frac{\beta+1}{(x_i+1)^2}+
\frac{\alpha+1}{(x_i-1)^2}  + 2 \sum_{1\leq j \leq n;\, j\neq i}
   \frac{1}{(x_i - x_j)^2}}& \text{if }i=j.
  \end{cases}
$$
It is easy to observe that matrix $H$ is real, symmetric, strictly
diagonally dominant, and its diagonal entries are positive, from
which we conclude that $H$ is positive definite. In particular,
every critical point of $E$ must be a point of a local minimum. Its
uniqueness was observed above, hence it is in fact the point of the
global minimum.

Stieltjes considered also similar electrostatic models for other two
classical orthogonal polynomials: Laguerre and Hermite. Since in
this situation the free charges can move on an unbounded set, what
can't prevent them from escaping to infinity? Stieltjes found a
clever solution in putting a constraint on either the first
(Laguerre) or second (Hermite) moment of their zero counting
measures (see e.g.\ \cite{Walter93}).

Namely, if for $\alpha>-1$ we fix one charge of mass $(\alpha+1)/2$
at the origin, and allow $n$ positive unit charges $X$ to move in
$[0,+\infty)$ with an additional constraint that their arithmetic
mean is uniformly bounded, $\sum_{i=1}^n x_i/n\leq K$ for some
positive constant $K$, then again the unique configuration providing
the global minimum of the total energy of the system of interacting
particles coincides with the rescaled zeros of the Laguerre
polynomial $L^{(\alpha)}_n(r_n x)$, where $r_n=(n+\alpha)/K$. In a
similar way, allowing $X$ to move freely on the whole $\R$, but
limiting the arithmetic mean of their \emph{squares}, $\sum_{i=1}^n
x_i^2/n\leq K$ for some constant $K$, the unique configuration
providing the global minimum of the total energy of the system of
interacting particles coincides with the rescaled zeros of the
Hermite polynomial $H_n(s_nx)$, where $s_n=\sqrt{(n-1)/(2K)}$.

The proof for these both cases is similar to the Jacobi case, except
that now we are dealing with a constrained minimum, so that the
Lagrange multipliers corresponding to the restrictions when looking
for a necessary condition of a critical point will become part of
the characterizing differential equation.

It is strange enough that Stieltjes himself did not realize that the
unbounded cases (Laguerre and Hermite) could have been put in the
same framework that the Jacobi one if we allow for an external field
not necessarily generated by positive fixed charges. It was probably
Ismail \cite{Ismail2000} who observed first that the zeros of these
polynomials still provide the global minimum of the total energy
\eqref{energia} without any additional constraint or rescaling, if
assuming for the Laguerre polynomials that $\varphi$ is generated by
a fixed charge of mass $(\alpha+1)/2$ at the origin \emph{plus}
another term, equal to $x/2$; analogously, for the Hermite
polynomials it is sufficient to take $\varphi(x)=x^2/2$.

Reviewing the electrostatic models above several natural questions
arise, such as:
\begin{itemize}
  \item are there generalizations of these models to other families
  of polynomials?
  \item why necessarily the global minimum of the energy should be
  considered? Which other types of equilibria described above could
  be linked to the zeros of the polynomials?
  \item what is the appropriate model for the complex zeros (when
  they exist)?
  \item what kind of applications can these models have beyond their
  clear aesthetical value?
\end{itemize}
In the following sections we will try to give at least some partial
answers or to formulate conjectures related to these questions.

%%%%%%%%%%%%%%%%%%%%%%%%%%%%%%%%%%%%%%
\section{A generalization of the electrostatic model: Lam\'{e} equation} \label{sec3}
%%%%%%%%%%%%%%%%%%%%%%%%%%%%%%%%%%%%%%

In fact, Stieltjes studied a more general situation. The
\emph{generalized Lam\'{e} differential} equation (in an algebraic form)
is
\begin{equation}\label{lame1}
   E''(x)+\left( \sum_{i=0}^p \frac{\rho_i}{x-a_i} \right) E'(x)-\frac{C(x)}{A(x)}\,
   E(x)=0\,, \quad A(x)=\prod_{i=0}^p (x-a_i)\,,
\end{equation}
where $C$ is a polynomial of degree $\leq p-1$ (in the sequel we use
the notation $C \in \P_{p-1}$). The case $p=1$ corresponds to the
hypergeometric differential equation, such as (\ref{EDJac2}), whyle
for $p=2$ we obtain the Heun's equation (see \cite{Ronveaux95}).

Heine \cite{Heine1878} proved that for every $n \in \N$ there exist
at most
\begin{equation}\label{combinatorial}
\sigma(n)=\binom{n+p-1}{n}
\end{equation}
different polynomials $C$ in (\ref{lame1}) such that this equation
has a polynomial solution of degree $n$. These coefficients $C$ are
called \emph{Van Vleck polynomials} and the corresponding polynomial
solutions $E$ are known as \emph{Heine-Stieltjes polynomials}.

Stieltjes studied the problem under the following two assumptions,
generalizing the classical situation for the Jacobi polynomials:
\emph{(i)} the zeros $a_i$ of $A$ are assumed to be simple and real,
so without loss of generality we can take
\begin{equation}\label{simplezeros}
-1=a_0 < a_1< \dots < a_p=1\,,
\end{equation}
\emph{(ii)} the residues of the coefficients of $E'$ are all
positive:
\begin{equation}\label{lame2}
 \rho_i>0, \quad i=0, \dots, p\,.
\end{equation}
(cf.\ equation (\ref{EDJac2}) with assumptions $\alpha >-1$, $\beta
>-1$). Observe that zeros $a_k$'s define $p$ intervals
$(a_{k-1},a_{k})$ of the real line, and the combinatorial number
$\sigma(n)$ defined in \eqref{combinatorial} coincides with the
number of different ways of placing $n$ balls in $p$ boxes. This is
not a mere coincidence. According to Stieltjes \cite{Stieltjes:885},
for each vector $\vt{n}=(n_1, \dots, n_p)$ such that
$n_1+\dots+n_p=n$ there exists a unique Heine-Stieltjes polynomial
$E_{\vt{n}}$ and the unique corresponding Van Vleck polynomial
$C_{\vt{n}} \in \P_{p-1}$ such that precisely $n_k$ zeros of
$E_{\vt{n}}$ belong to the interval $(a_{k-1},a_{k})$, $k=1, \dots,
p$. Moreover, these zeros provide the \emph{unique global minimum}
for the total energy \eqref{energia} under the \emph{additional
restriction} of the number of zeros in each interval mentioned
above, with the external field created by $p+1$ positive charges
fixed at $a_j$'s:
\begin{equation}\label{externali_field_Heine}
\varphi(x)=-\sum_{j=0}^p \frac{\rho_j}{2} \, \ln|x-a_j|\,.
\end{equation}

This electrostatic model has been used in \cite{MS:01} in order to
study the semiclassical limit $n\to \infty$ (see Section
\ref{sec5}), and in a number of papers \cite{BT:01,BJMT:03} in
describing the probabilistic distribution of these zeros as the
number of intervals (that is, $p$) grows large together with $n$.

Further generalizations of the work by Heine and Stieltjes followed
several paths; we will mention only some of them. First, under
assumptions (\ref{simplezeros})--(\ref{lame2}) Van Vleck
\cite{Vleck1898} and B\^{o}cher \cite{Bocher97} proved that the zeros of
$C$ belong to $[a_0,a_p]$. A refinement of this result is due to
several works of Shah \cite{Shah:68,Shah:69b,Shah:69,Shah:70}.
Furthermore, P\'{o}lya \cite{Polya12} showed that, allowing complex
zeros of $A$, condition (\ref{lame2}) is sufficient to assure that
the zeros of $E$ are located in the convex hull of $a_k$'s. Marden
\cite{Marden66} and later on, Alam, Al-Rashed, and Zaheer (see
\cite{Alam79}, \cite{Al-Rashed85}, \cite{Zaheer76}, \cite{Zaheer77})
established further results on location of the zeros of the
Heine-Stieltjes polynomials under weaker conditions on the
coefficients $A$ and $\rho_i$ of (\ref{lame1}).

\begin{figure}[hbt]
\centering
\begin{tabular}{c@{\quad}c}
\mbox{\includegraphics[scale=0.6]{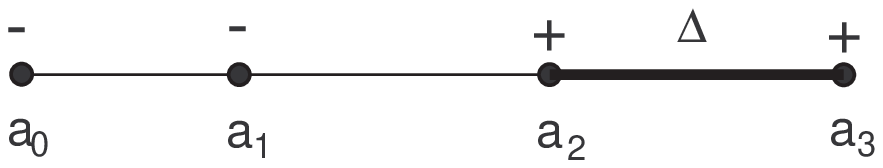}} &
\mbox{\includegraphics[scale=0.6]{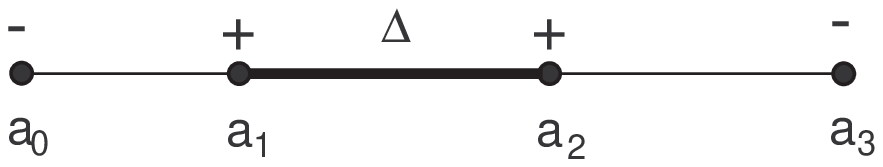}}
\end{tabular}
\caption{Configuration studied by Dimitrov and Van Assche.}
\label{fig:cases}
\end{figure}

Going back to configuration \eqref{simplezeros}, not much can be
said if we drop the condition of positivity of the residues $\rho_i$
in (\ref{lame2}); the situation then becomes much more difficult to
handle. The first attempt was made by Gr\"{u}nbaum
\cite{Grunbaum98,Grunbaum01}, and later on by Dimitrov and Van
Assche \cite{Dimitrov:00}. In particular, in \cite{Dimitrov:00} the
case of $p=3$ is analyzed when the positive and negative residues do
not interlace: basically, the two cases considered are depicted in
Figure \ref{fig:cases}.
%, that is,
%$$
%  \text{ {either} }  \quad \rho_{0}, \rho_{1}< 0, \rho_{2}, \rho_{3}> 0  \quad \text{ {or} }  \quad
%\rho_{0}, \rho_{3} < 0,  \rho_{1}, \rho_{2} > 0 .
%$$
%%%

%%%%%%%%%%%%%%%%%%%%%%%%%%%%%%%%%%%%%%%%%%%%%%%%%%%%%%%%%%
%%%%%%%%%%%%%%%%%%%%%%%%%%%%%%%%%%%%%%%%%%%%%%%%%%%%%%%%%

Dimitrov and Van Assche proved that under certain restrictions on
the degree of the polynomial and the residues $\rho _i$, for every
$n \in \N$ there exists a unique pair $(C_n,E_n)$  of respectively
Van Vleck and Heine-Stieltjes polynomials  with $\deg E_n=n$. All
zeros of $E_n$ belong to the interval enclosed by $a_j$'s with
$\rho_j>0$ (denoted by $\Delta $ in Figure \ref{fig:cases}), and
they are in the equilibrium position, given by the global minimum of
the discrete energy \eqref{energia} with the external field
\eqref{externali_field_Heine}. As far as we know, no further studies
in this direction exist.

Interestingly enough, a similar situation appears when we study a
special class of polynomials of multiple orthogonality (when the
orthogonality conditions are distributed among several measures),
which are object of an intensive study in the last years. Kaliaguine
\cite{kaliaguine:1981} (see also the more recent papers by Aptekarev
\cite{Aptekarev:97} and Kaliaguine and Ronveaux
\cite{kaliaguine:1996}) considered the system of polynomials $\{
p_{2n}\}$, $\deg p_{2n} =2n$, satisfying
$$
\int_{-1}^0 p_{2n}(x) x^k w(x)\, dx=\int_{0}^1 p_{2n}(x) x^k w(x) \,
dx=0\,, \quad k=0, 1, \dots, n-1\,,
$$
with $w(x)=x^\gamma (1-x)^\alpha  (1+x)^\beta $, $\alpha , \beta ,
\gamma >-1$. It is a real case: $p_{2n}$ has exactly $n$ zeros in
each interval $(-1, 0)$ and $(0,1)$. The asymptotic analysis (using
the tools of the potential theory) suggests that a reasonable model
for these zeros might be that each of these sets of $n$ zeros is the
equilibrium position of the positive unit charges, which besides the
mutual interaction are subject to the external field generated by
the three fixed charges $\alpha /2$, $\beta /2$, $\gamma /2$ at $1,
-1, 0$, respectively, and by the $n$ charges of the weight $1/2$ on
the other subinterval. Nevertheless, this is an open problem.

In a recent paper Grinshpan \cite{Grinshpan:01} considered a
configuration of zeros on the unit circle. He characterized the
stable equilibrium of the $n$ positive unit charges on $\T=\{z \in
\C:\, |z|=1 \}$ in the field generated by $n$ negative unit charges
in $\C \setminus (\T \cap \{0\})$ and showed that they are zeros of
a Heine-Stieltjes polynomial satisfying a generalized Lam\'{e}
differential equation  with coefficients symmetric with respect to
$\T$. Although it might give an impression that here we are dealing
with complex zeros, the symmetry of the problem makes it essentially
``real''. A more recent paper \cite{MF:05} extends the problem of
stable equilibrium studied in \cite{Grinshpan:01} to the case of a
field generated by $m=m(n)\in \N$ negative charges in $\C \setminus
(\T \cap \{0\})$ with values $-\omega_{nk}$, where $\omega_{nk}>0$
are the residues of $R_n/Q_n$ in the zeros $z_{nk}$ of $Q_n$ that
lie in the open unit disc, and $Q_n$, $R_n$ are polynomials of
degree $2m$, $2m-1$, respectively, such that, $A_n(z)=zQ_n(z)$ and
$B_n=Q_n(z)-zR_n(z)$ are the coefficients of the corresponding
generalized Lam\'{e} differential equation $A_n(z)y''(z)+B_n(z)y'(z)+
C_n(z) y(z)=0$.
%There asymptotic properties, where $N \to \infty$, of the
%Heine-Stieltjes and Van Vleck polynomials are studied.

%%%%%%%%%%%%%%%%%%%%%%%%%%%%%%%%%%%%%%
\section{Electrostatics and orthogonality} \label{secOrthogonality}
%%%%%%%%%%%%%%%%%%%%%%%%%%%%%%%%%%%%%%

Stieltjes' theorem is beautiful, and also useful for proving some
monotonicity properties of the zeros with respect to the parameters
of the system. But is the existence of such a description a mere
accident? Why should an electrostatic interpretation of the zeros of
general orthogonal polynomials exist at all? One fact (known from
the '80) that makes us expect further developments is that the
appropriate description of the asymptotic behavior (as $n\to
\infty$) of the zeros of orthogonal polynomials is given in terms of
some \emph{equilibrium measures}, that is, measures minimizing
certain logarithmic potentials, eventually in presence of an
external field (see \cite{Saff:97}).

If we go back to the Ismail's electrostatic model for the Laguerre
and Hermite polynomials, then it is easy to recognize in the
corresponding external field the influence of the weight of
orthogonality of these polynomials. It was definitely an interesting
problem to find out if this is a regular fact, tackled by Ismail in
a series of papers \cite{Ismail2000,Ismail2000_2}. We have noticed
that so far the characterizing differential equation has been the
cornerstone of the model, and it was a reasonable line of attack.
Thus, generalizing some pioneering works of Bauldry
\cite{Bauldry:90}, and Bonan and Clark \cite{Bonan:90}, Ismail
deduced a second order linear differential equation satisfied by a
sequence of \emph{general} orthogonal polynomials, after which he
suggested a model for the zeros. Let us outline briefly his results.

Let $\{p_n(x)\}$ be polynomials orthonormal with respect to a unit
weight function $w(x)=\exp(-v(x))$ supported on an interval $[c , d
]\subset \R$, finite or infinite:
$$
\int_c ^d  p_m(x) p_n(x) w(x) \, dx = \delta_{m,n}\,.
$$
Then sequence $\{p_n\}$ satisfies a three-term recurrence relation
$$
x p_n(x) =a_{n+1}p_{n+1}(x) +b_n p_n(x) +a_n p_{n - 1}(x), \quad n
\geq 0, \quad a_n>0,
$$
with $p_{-1}(x) = 0$ and $ p_0(x) = 1$. Moreover, it has been shown
in \cite{Bauldry:90}, \cite{Bonan:90}, and also \cite{Chen:97},
\cite{Ismail98}, that under certain assumptions on $w$ the
orthonormal polynomials $p_n$ satisfy also the
difference-differential relation
$$
 p'_n (x) =A_n(x)p_{n - 1}(x) - B_n(x)p_n(x)\,,
 $$
where the coefficients $A_n$ and $B_n$ are explicitly given in terms
of $w$, the recurrence coefficients $a_n$, and the values of $p_n$
at the endpoints $c $ and $d $. A direct consequence of this is that
$p_n$ satisfies also the second order linear differential equation
\begin{equation}\label{diff_eq}
p''_n (x) - 2 R_n(x)p'_n (x) +S_n(x)p_n(x) = 0\,,
\end{equation}
with
\begin{align*}
R_n(x) &=   \frac{v'(x)}{2} + \frac{A'_n(x)}{2A_n(x)}  \,, \\
S_n(x) &=B'_n (x) - B_n(x) \frac{A'_n (x)}{ A_n(x)} - B_n(x)[v'(x)
+B_n(x)] + \frac{a_n}{ a_{n-1}}\,  A_n(x) A_{n - 1}(x)\,.
\end{align*}
In particular, if $X^*\subset (\alpha ,\beta )^n$ is the ordered set
of zeros of $p_n$, then
\begin{equation}\label{necessary_Ismail}
p''_n (x) - 2 R_n(x)p'_n (x) = 0  \quad \text{for }   x \in X^*\,.
\end{equation}
Comparing it with \eqref{equations_critical2} can associate
$\varphi'(x)=R_n(X)$, and it is natural to consider an electrostatic
model for $X^*$ with logarithmic interaction between particles and
an external field
\begin{equation}\label{external_Ismail}
\varphi(x)= \frac{v(x)}{2} + \frac{\ln \left(k_n
A_n(x)\right)}{2}=\varphi_{\rm long}(x)+\varphi_{\rm short}(x)
\end{equation}
($k_n$ is any appropriate normalization constant, taken in
\cite{Ismail2000} equal to $a_n^{-1}$). Observe that this external
field has two components: the first term in the right hand side of
\eqref{external_Ismail} has its origin in the orthogonality weight
$w(x)=\exp(-v(x))$, and Ismail called it the \emph{long range
potential}. The second term, which received the name of the
\emph{short range potential}, is a bit mysterious, but explains
several features of the classical models considered in Section
\ref{sec2}, and at the same time, allows to give a further
generalization of the electrostatic interpretation.

The main result of Ismail \cite{Ismail2000} states that, assuming
$w(x) > 0$ on $ (c , d )$, both $v$ and $\ln(A_n)$ in $C^2(c ,d )$,
and the external field \eqref{external_Ismail} convex, the total
energy \eqref{energia} has a unique point of global minimum, which
is precisely $X^*$ (that is, the zeros of the orthogonal polynomial
$p_n$). Obviously, what is used here is in fact the relation
\eqref{necessary_Ismail} equivalent to the fact that $X^*$ is just a
critical point for $E(X)$, and convexity assures that this is the
unique point of minimum. The proof is in general very similar to the
arguments for the Jacobi polynomials from Section \ref{sec2}.

Beside some restrictions on the weight of orthogonality, one
drawback of this model in application to general weights is that
usually the short range potential cannot be explicitly computed in
terms of $w$. This is not the case of the classical polynomials and
some of their generalizations.

For instance, for Jacobi polynomials,
  $$
  w(x)=\frac{(1-x)^\alpha (1+x)^\beta \Gamma (\alpha +\beta +2)}%
  {2^{\alpha +\beta +1}\Gamma (\alpha +1)\Gamma (\beta +1)}\,,  \quad
  \frac{A_n(x)}{a_n}=\frac{2n+\alpha +\beta +1}{1-x^2}\,.
  $$
In other words, $w$ is responsible for the fields generated by two
positive charges of mass $\beta/2$ and $\alpha/2$ at $-1$ and $+1$,
respectively, while the short range potential adds the missing
charges of size $1/2$ at $\pm 1$.

Analogously, for Laguerre polynomials we have
  $$
  w(x)=\frac{x^\alpha \exp(-x)}{\Gamma (\alpha +1)}\,,  \quad
  \frac{A_n(x)}{a_n}=\frac{1}{x}\,,
  $$
so that the $x^\alpha $ factor in $w$ plus the short range potential
are responsible for the fixed charge of size $(\alpha +1)/2$ at
$z=0$, and the exponential factor of the weight $w$ generates the
remaining component $x/2$ of the external field $\varphi$, as
described at the end of Section \ref{sec2}.

So far the only role played by the short range component was adding
some extra weight to the charges fixed by the weight. However, $A_n$
can be responsible for a creation of ``ghost'' movable charges
(either positive or negative) in the picture, as it was shown first
by Gr\"{u}nbaum \cite{Grunbaum98} (work that in fact predates
\cite{Ismail2000}) or by Ismail himself when analyzing the case of
the so-called Freud weight, $w(x)=\exp(-x^4)$.

Let us consider here a little bit more general situation. A
generalized weight function (or measure) is \emph{semiclassical}
(see e.g.\ \cite{Maroni91,Shohat:39}) if it satisfies the
\emph{Pearson equation}
\begin{equation*}%\label{pearson}
D(\phi w)=\psi w\,,
\end{equation*}
where $\phi $, $\psi$ are polynomials, with degree of $\psi \geq 1$,
and $D$ is the ``derivative'' operator (in  the usual, but also
possibly in a distributional sense). It is  well known that for such
a weight the corresponding orthogonal polynomials (called also
\emph{semiclassical}) satisfy a differential equation of the type
(\ref{diff_eq}), where the coefficients $R_n$ and $S_n$ are rational
functions. The classical-type orthogonal polynomials considered by
Gr\"{u}nbaum \cite{Grunbaum98} are an example of a semiclassical family,
but there are many more (for instance, the so-called sieved
orthogonal polynomials \cite{Ismail98}, among others).

The systematic study of the semiclassical orthogonal polynomials
from the point of view of the electrostatic interpretation of their
zeros has started in the works of Garrido and collaborators
\cite{Garrido05}, \cite{Garrido_Tesis}. For instance, they consider
a perturbation of the Freud weight function ($w(x)=\exp(-x^4)$) by
the addition of a fixed charged point of mass $\lambda$ at the
origin; the corresponding orthogonal polynomials are known as
Freud-type polynomials. The resulting orthogonality measure is
semiclassical, and it was proved in \cite{Garrido05} that these
polynomials satisfy a second order linear differential equation of
the form \eqref{diff_eq}, and the electrostatic model is in sight.
In fact, in the situation considered, and following Ismail's
terminology, the long range potential is, as expected, $\varphi_{\rm
long}(x)  = x^4/2$, while the short range potential depends on the
degree $n$ and has the following structure:
\begin{equation}\label{short_Garrido}
\varphi_{\rm short}(x)= \frac{1}{2}\, \ln\left|\frac{
 Q_4(x,n) }{x^2}\right|
\end{equation}
where $Q_4(x,n)$ has two real roots $r_1(n), r_2(n)$ and two simple
conjugate complex roots $r_3(n), r_4(n)$. Their asymptotic behavior
is
\begin{equation}\label{limits_charges}
\lim_{n} r_1(n) = \lim_{n} r_2(n)=0\,, \quad \lim_{n} r_3(n) =
\lim_{n} r_4(n)=\infty\,.
\end{equation}
In fact, in the odd case, both real zeros coalesce at the origin:
$$
r_1(2n+1)=r_2(2n+1)=0\,,
$$
canceling out the denominator in the argument of the logarithm in
\eqref{short_Garrido}. In other words, the short range potential is
given by one unit positive charge fixed at the origin, plus 4
floating negative charges of weight $-1/2$ on $\C$ situated
symmetrically with respect to the origin. However, in the odd case
two of these negative charges annihilate with the fixed positive
one, paving the way to a zero of the orthogonal polynomial to take
the origin, position it must occupy due to the symmetry of the
problem. This annihilation has as a consequence that for odd $n$'s
the Freud and the Freud-type polynomials are identical. Moreover,
\eqref{limits_charges} shows that also for the even case they are
asymptotically identical (as $n\to \infty$).

One of the main results of \cite{Garrido05}, \cite{Garrido_Tesis} is
that the zeros of the Freud-type polynomials provide a
\emph{critical configuration} for the total energy $E(X)$ (in other
words, that for these zeros $X^*$ equations
\eqref{equations_critical} are satisfied). But can we assure that
they are in a stable equilibrium, that is, that  $E(X^*)=\min E(X)$?
The complete answer is not clear yet. By computing the Hessian
matrix $H$ \eqref{hesssian} it was proved that for small values of
the correction charge $\lambda $ this is really the case, but there
are some evidences that make it possible to conjecture the existence
of a critical value $\lambda _0$ such that for
$0\leq\lambda<\lambda_0(n)$ the zeros of the orthogonal polynomials
attain the minimum of the total energy, but for $\lambda
>\lambda _0(n)$ the equilibrium is no longer stable. Is still
any other type of equilibria discussed in Section \ref{sec2}
preserved in this case (beside the critical point of the energy)?
This is also an open question.

%%%%%%%%%%%%%%%%%%%%%%%%%%%%%%%%%%%%%%
\section{A max-min problem} \label{secMaxMin}
%%%%%%%%%%%%%%%%%%%%%%%%%%%%%%%%%%%%%%

From potential theory it is well known that many of energy
minimization problems have in fact a min-max (or max-min)
characterization. Let us analyze again, but from this point of view,
the classical models of Stieltjes, and in passing, address the issue
why the real line should be present in the model.

This question is not trivial: let us go back to the Jacobi
polynomials $P_n^{(\alpha,\beta)}$ with $\alpha , \beta
>-1$. Why their zeros should be real? A standard explanation is that
these polynomials satisfy the orthogonality conditions
\eqref{orthogonalityForJacobi}, but where the real line is present
in that integral? Observe that the integrand in
\eqref{orthogonalityForJacobi} is an analytic function, so we
perfectly can obtain the same result integrating along any
reasonable curve joining $-1$ with $+1$.

The real line is also built in the electrostatic model of Stieltjes
as an a priori constraint. However, we may put forward an
alternative model free from this restriction. As in the Stieltjes'
setting we fix two positive charges of mass $(\beta+1)/2$ and
$(\alpha+1)/2$ at $-1$ and $+1$, respectively, and denote again by
$E(X)$ the total energy of a discrete system $X\subset \C^n$, given
by formulas \eqref{external_field_Jacobi}--\eqref{energia} (assuming
that in case of coincidence of any two points, $E(X)=+\infty$).
Denote by $\F$ the family of compact continua on $\C$ containing
both $- 1$ and $+1$.  For any $K\in \F$ set
$$
m(K)=\inf\{E(X):\, X=\{x_1, \dots, x_N \}\subset K \}\,,
$$
and consider the following extremal problem:
\begin{equation*}%\label{problemaExtremal}
E^*=\sup\{m(K):\, K \in \F \}\,.
\end{equation*}
\begin{thm} \label{theorem1}
$E^*<+\infty$, and this value is attained at a unique $n$-tuple
$X^*=\{x_1^*,\dots, x_n^* \} $, characterized by the fact that
$x_i^*$ are precisely the zeros of the Jacobi polynomial
$P_n^{(\alpha , \beta )}$.
\end{thm}
\begin{proof}
Let $X$ be any configuration of $n$ distinct points on the interval
$(-1,1)$ (the case of coalescence of some points is trivial), and
let $K\subset \F$ be any other continuum. Since $\pm 1 \in K$, we
can assure that there exists a discrete set of points $Z\subset K$
of the form
$$
Z=\left \{z_j=x_j+i y_j:\, y_j\in \R, \; j=1, \dots, n \right\}
$$
(that is, any vertical line passing through an $x_j$ must intersect
$K$). It is straightforward to check that
$$
|x_k-x_j|\leq |z_k-z_j|\,, \quad k\neq j\,, \quad \text{and} \quad
|1\pm x_j|\leq |1\pm z_j|\,, \quad j=1, \dots, n\,.
$$
Taking into account the expression of the total energy and the
positivity of all the charges involved in the system we immediately
conclude that $ E(X)\geq E(Z)$. Hence,
$$
m([-1,1]) \geq m(K)\,,
$$
and since $K$ was arbitrary, we obtain that $E^*=m([-1,1])$. But
from Stieltjes' theorem it follows that there exists a unique
$X^*\subset (-1,1)^n$ such that $E(X^*)=m([-1,1])$, and the points
of $X^*$ are the zeros of $P_N^{(\alpha , \beta )}$, which concludes
the proof.
\end{proof}

Obviously, a similar result is possible to formulate for Laguerre
and Hermite polynomials.

\begin{figure}[hbt]
\centering
\includegraphics[scale=0.6]{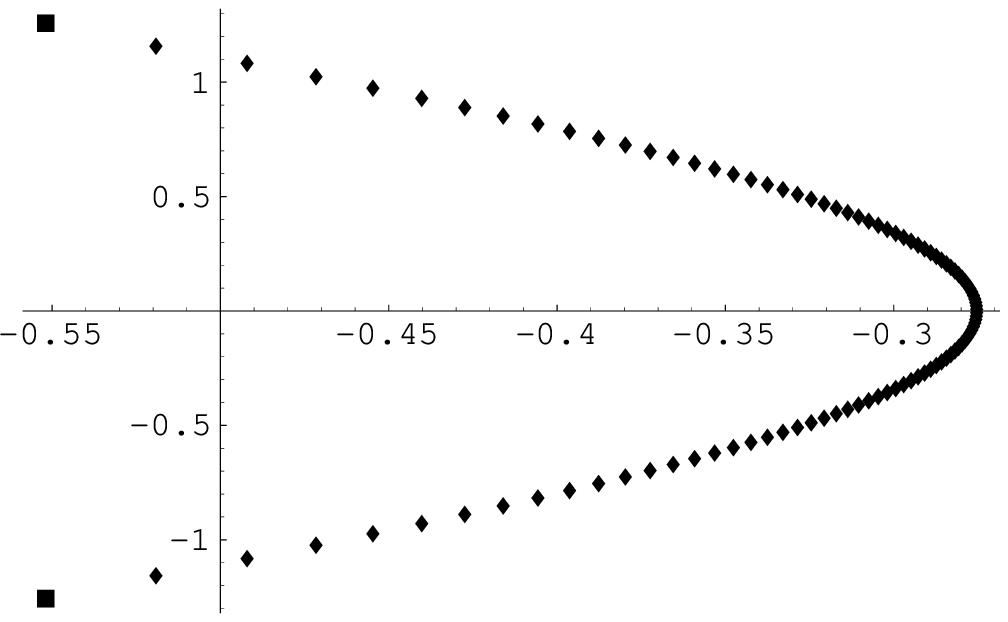}
\quad
\includegraphics[scale=0.6]{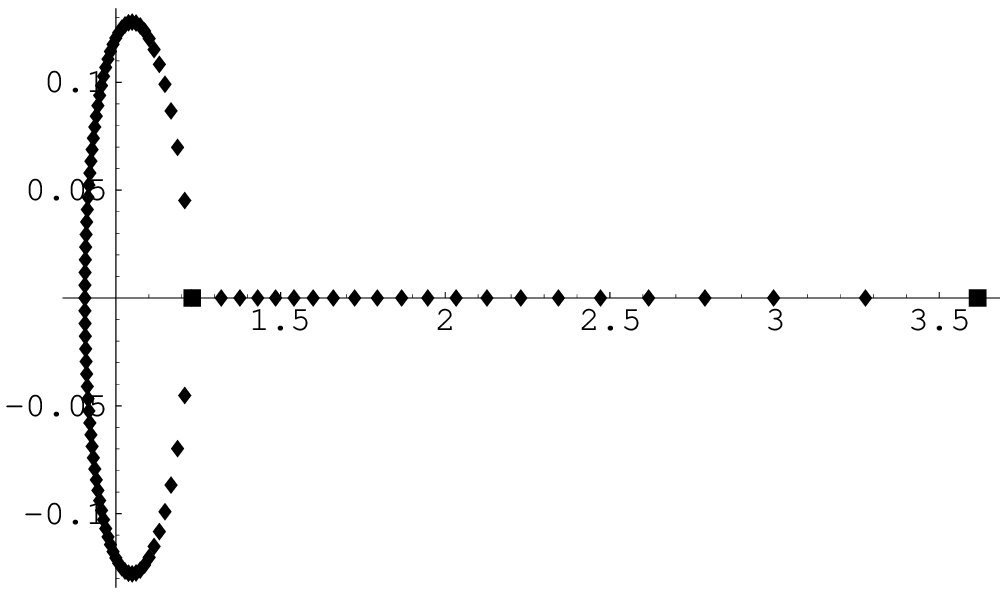}
\caption{\textbf Zeros of $P_{100}^{(-300,-200)}$ and
$P_{100}^{(-80,-314)}$. } \label{fig2}
\end{figure}

The value of the min-max model we have described is that it allows
to consider easily complex zeros of these classical families, which
appear when the parameters involved become no longer ``classical''
(see for instance Figure \ref{fig2} with two examples of zeros of
Jacobi polynomials). In fact, as it follows from the Rodrigues
formula \eqref{rodrigues}, if we allow for arbitrary real $\alpha$
and $\beta $, the Jacobi polynomials constitute a very rich object
exhibiting several types of orthogonalities: hermitian,
non-hermitian, and even multiple, depending on the values of the
parameters and the degree (see the classification of all cases
obtained \cite{KMF:05}). We conjecture that Theorem \ref{theorem1}
is valid in all these situations after a suitable modification of
the family of continua $\F$. For instance, if $2n+\alpha +\beta <-1$
(see Figure \ref{fig2}, left), the appropriate class $\F$ would
include all continua extending to infinity, and \emph{separating}
$-1$ and $+1$. From the results of \cite{KMF:05} a rule of the thumb
follows: if $\alpha +n>0$ (resp., if $\beta +n>0$) then the continua
from $\F$ must contain $+1$ (resp., $-1$). Moreover, if $n+\alpha
+\beta +1<0$, then the continua from $\F$ must include the infinity.

%%%%%%%%%%%%%%%%%%%%%%%%%%%%%%%%%%%%%%%%%%%%%%%%%%%%%%%%%%
%%%%%%%%%%%%%%%%%%%%%%%%%%%%%%%%%%%%%%%%%%%%%%%%%%%%%%%%%%

All these facts have been conjectured by one of this authors some
time ago, but their proof still remains an open problem. A general
approach to such kind of theorems would allow also to tackle other
kind of problems not linked directly to the real line or exhibiting
non-real zeros.

Typical examples are those related to the multiple (or Hermite-Pad\'{e})
orthogonality, when the orthogonality conditions are distributed
among several measures, or of non-hermitian orthogonality, when the
path of integration is, in fact, not fixed a priori. In many
situations these polynomials exhibit complex zeros, whose
\emph{asymptotic} distribution has been studied by several authors.
From the pioneering works of Stahl \cite{Stahl:86} it is known that
the limiting location of the zeros is usually described by some
trajectories of a quadratic differential, which in turn have a
min-max description. This gives an additional evidence to our
conjecture, that in fact says that this behavior is not exclusive of
large $n$'s (a kind of ``discrete version'' of Stahl's theorem).

%%%%%%%%%%%%%%%%%%%%%%%%%%%%%%%%%%%%%%
\section{Further generalizations and applications
of the electrostatic model} \label{secFurther}
%%%%%%%%%%%%%%%%%%%%%%%%%%%%%%%%%%%%%%

Many new interesting problems related to the electrostatic models
have their origin either in some non-standard orthogonality
conditions (such as non-hermitian or multiple orthogonality
mentioned above), or in physical applications. In this Section we
will briefly review some of them.

Exactly solvable or quasi-exactly solvable multi-particle quantum
mechanical systems have many remarkable properties. Especially,
those of the Calogero-Sutherland-Moser (CSM) systems
\cite{Calo71,Calo75,Moser75,Suth72} and their integrable deformation
called the Ruijsenaars-Schneider-van Diejen (RSvD) systems
\cite{RS86,vD94} have been well studied. A classical result is that
the equilibrium positions of the CSM systems are described by the
zeros of the classical orthogonal polynomials; the Hermite,
Laguerre, and Jacobi polynomials \cite{Calo77a,Calo77b} (see also
\cite{OS:05} for a comprehensive review and bibliography therein).
Following this analogy the authors of \cite{OS:05} proved recently
that the equilibrium positions of the RSvD systems with rational and
trigonometric potentials coincide again with zeros of polynomials
from the Askey tableau of hypergeometric orthogonal polynomials;
namely, they found connections with the Meixner-Pollaczek,
continuous Hahn, continuous dual Hahn, and the Askey-Wilson
polynomials. Hence, an electrostatic model for the zeros of
hypergeometric polynomials in all their generality could have
several interesting applications for exactly solvable
quantum-mechanical systems.

Another interesting model was considered by Loutsenko
\cite{loutsenko:2003}, who studied a system $(X,Y)\subset \C^{n+m}$
of $n$ positive and $m$ negative \emph{moving} charges in $\C$ of
masses 1 and $-\Lambda$, respectively (with $\Lambda\in \{
1/2,1,2\}$), interacting again according to the logarithmic law. He
proved in particular that if $(X^*,Y^*)$ is a critical configuration
of the mutual energy $E_{\rm mutual}(X,Y)$, $p$ is the monic
polynomial vanishing at the points of $X^*$, and $q$ is the monic
polynomial vanishing at those of $Y^*$, then $p, q$ are solutions of
a \emph{bilinear} differential equation $$p'' q-2\Lambda p'q'+
\Lambda^2 p\, q''=0.$$

Regarding the Lam\'{e} equation, once we depart from the classical
setting considered by Stieltjes, the description of the
Heine-Stieltjes polynomials is not apparent. It is clear however
that we can hardly expect that their zeros provide a stable
equilibrium; instead, we must concentrate on the least restrictive
condition of configuration providing critical values of the energy
functional. This study has again some physical applications. The
$BC_n$ elliptic Inozemtsev model is a quantum integrable system
with $n$-particles whose potential is given by elliptic functions.
For the case $n=1$, finding eigenstates of its Hamiltonian is
equivalent to solving the Heun equation (see e.g.\
\cite{Tak00,Tak03}. In this sense, the $BC_n$ Inozemtsev model is a
generalization of the Heun equation.

There are also connections of some problems in physics and
representation theory with multiple orthogonal polynomials. It is
known that zeros of the Jacobi polynomial $P_n^{(\alpha , \beta )}$
satisfy a system of algebraic equations, which is known as the
\emph{Bethe Ansatz equation} of the Gaudin model associated to
$sl_2$  and two irreducible modules with highest weights $-(\alpha
+1),-(\beta+1)$. In \cite{MV:05b} the authors generalized this
connection, studying sequences of $r$ polynomials whose zeros
constitute the unique solution of the Bethe Ansatz equation
associated with two highest weight $sl_{r+1}$ irreducible modules,
with the restriction that the highest weight of one of the modules
is a multiple of the first fundamental weight. As a result, they
show that the first polynomial in the sequence coincides with the
well known Jacobi-Pi\~{n}eiro multiple orthogonal polynomial, and others
are given by Wronskian type determinants of Jacobi-Pi\~{n}eiro
polynomials. In \cite{MV:05a} they derived a linear differential
equation (alas, of order $r$) for these polynomials, paving a way to
the electrostatic interpretation of their zeros.

Finally, we can mention the matrix orthogonal polynomials, whose
non-trivial connections with (matrix) differential equations is
being studied. Duran and Gr\"{u}nbaum  gave  in \cite{Duran:04,Duran:05}
some  examples of matrix orthogonal polynomials $\{P_n\}$ satisfying
the second order linear ODE of the form
$$
P''_n(x)A_2(x)+P'_n(x)A_1(x)+P_n(x)A_0(x)=\Gamma_n P_n(x),\quad
n\geq 0\,,
$$
where (as in the scalar hypergeometric case) coefficients $A_j$ are
matrix polynomials that do not depend on $n$, of degrees $\leq 2$,
$1$, and $0$, respectively, and $\Gamma_n$ are Hermitian matrices.
These $\{P_n\}$ look like a natural generalization of the families
of classical polynomials in the scalar case; likewise,  they exhibit
a rich variety of structural properties (see \cite{Duran:toappear}).
Furthermore, they might have some relevance in the analysis of the
Dirac equation (relativistic analogue of the Schrodinger equation).
However, a study of the electrostatic interpretation of the zeros of
the orthogonal matrix polynomials remains completely open.

%%%%%%%%%%%%%%%%%%%%%%%%%%%%%%%%%%%%%%
\section{Asymptotic distribution of zeros} \label{sec5}
%%%%%%%%%%%%%%%%%%%%%%%%%%%%%%%%%%%%%%

Assume we have an electrostatic model of $n$ positive charges moving
in an external field, and we are interested to analyze what happens
when $n\to \infty$ (the so-called thermodynamic or semiclassical
limit). In many situations this information can be extracted
directly from the characterizing second order linear differential
equation, using either the WKB approach (cf.\
\cite{Arriola:85}--\cite{Arriola:91}, \cite{MF1:95}--\cite{PMG:01},
\cite{ZA:93}--\cite{ZA2:95}), or reducing the ODE to a Riccati form,
a method that we outline below. But first we should introduce some
notation and agree in the meaning of the ``global'' asymptotics of
the zeros.

In the sequel, $\supp(\mu)$ denotes the support of a measure $\mu$,
and
$$
\hat{\mu}(z)=\int\frac{d\mu(t)}{z-t}, \quad z\in \C\backslash
\supp(\mu)\,,
$$ is its Stieltjes (or Cauchy) transform. With a function $y: \C \rightarrow \C$ we can associate its
\textit{zero counting measure} $\nu_y$, defined by
$$
\nu_y = \sum_{y(x)=0} \delta_x\,,
$$
where $\delta_x$ is the unit mass (Dirac delta) at $x$, and the
sum goes along all the zeros of $y$ taking into account their
multiplicity. Equivalently, for each Borel set $\Delta \subseteq
\C$ the number of zeros of $y$ in $\Delta $ is
$$
\nu_y(\Delta) = \int_\Delta d\nu_y(x) \;.
$$
An important observation is that if $p$ is a polynomial and $\nu_p$
is the associated zero-counting measure, then the logarithmic
derivative of $p$ is the Stieltjes transform of $\nu$:
\begin{equation}\label{logDerivative}
\frac{p'(z)}{p(z)}=\sum_{k=1}^{n} \frac{1}{z-a_k}=\int
\frac{d\nu(t)}{z-t}=\hat{\nu}(z)\,,
\end{equation}
which can be evaluated away from the zeros of $p$.

The ``global'' behavior of the zeros of polynomials $\{ p_n\}$,
$\deg(p_n)=n$, is described by the limit of the sequence of
\emph{normalized} measures $\mu_n=\nu_{p_n}/n$ as $n \to \infty$, in
the sense of the weak-* topology. Recall (cf.\ \cite[Section
2]{Billingsley:71}) that a sequence of Borel measures $\{\mu_n\}$
converges to a measure $ \mu $ in the weak-* topology (which we
denote as $\mu _n \stackrel {*} {\longrightarrow} \mu$) if
 $$
 \lim_n \int f(x)d \mu_n (x) =  \int f(x)d \mu (x)\,, \quad \forall f
 \text{ continuous and compactly supported on } \C\,.
 $$
 Another important fact is that the set of unit measures with
 uniformly bounded supports is compact in the weak-* topology. Hence,
 if we know for instance that all the zeros of the sequence of
 polynomials $\{p_n\}$, $\deg(p_n)=n$, belong to the same compact
 set $K\subset \C$, then there always exists a unit measure $\mu$
 supported on $K$ and a subsequence $\Lambda\subset \N$ such that
 $\mu_n \stackrel {*} {\rightarrow}
 \mu$ for $n\in \Lambda$ (where $\mu_n$ are the normalized zero counting measures
 defined above). In consequence, and taking into account
 \eqref{logDerivative}, we have that
\begin{equation}\label{weakLimitJacobi}
\lim_{n \in \Lambda}\,  \frac{1}{n}\, \frac{p_n'(z)}{p_n(z)}=
\hat{\mu}(z)\,,\quad z\in \C \setminus \supp(\mu)\,.
\end{equation}
The method of reduction of the original second order ODE to the
Riccati form is based on the previous observation. This idea appears
in a work of Saff, Ullman, and Varga \cite{Saff:80}, although its
roots can be traced back to the famous Perron's monograph
\cite{Perron:50}, who in turn gives credit to some original works of
Euler. Recently it has been successfully applied in a variety of
problems (see, e.\ g.\ \cite{Faldey:95,MF2:99,MF1:99,PMG:01}).

Let us see how it works in the simplest case of Jacobi polynomials
$p_n=P_n^{(\alpha, \beta )}$ with $\alpha , \beta >-1$. The first
step is to rewrite equation (\ref{EDJac2}) in terms of the
normalized logarithmic derivative of the polynomial solution,
$$
h_n(z)=\frac{1}{n}\, \frac{p_n'(z)}{ p_n(z) }\,,
$$
which yields
\begin{equation}\label{beforeLimits}
\left(\frac{h'_n(z)}{n}+h_n^2(z)\right)+
\left(\frac{\alpha+1}{z-1}+\frac{\beta+1}{z+1}\right)\frac{h_n(x)}{n}+
\frac{n+\alpha+\beta+1}{n}=0\,.
\end{equation}
Since all the zeros of $p_n$ are in $[-1,1]$, by the argument
explained above we can assure the existence of a unit measure $\mu$
 supported on $[-1,1]$ and a subsequence $\Lambda\subset \N$ such
 that (see \eqref{weakLimitJacobi}) $h_n(z) \to \hat{\mu}(z)$, $n \in \Lambda $,
uniformly on compact subsets of $\C \setminus [-1,1]$, where $h_n'$
are also uniformly bounded. Taking limits in (\ref{beforeLimits})
when $n\in \Lambda$ we arrive at an \emph{algebraic} equation for
$h=\hat{\mu}$: $(1-x^2)h^2(x)+1=0$, so that
$$
\hat \mu(z)=(z^2-1)^{-1/2} \,,
$$
and the appropriate branch of the square root in $\C \setminus
[-1,1]$ is fixed by the condition $\lim_{z \to \infty} z \hat
\mu(z)=1$. Furthermore, since $\hat \mu$ is independent on
$\Lambda$, we can claim that $\mu$ is in fact the weak-* limit of
$\{\mu_n\}$ along $n\to \infty$.

It remains only to recover $\mu$ from its Stieltjes transform. Since
we know (important fact!) that $\supp(\mu)\subset [-1,1]$, this task
is just a straightforward application of the Sokhotsky-Plemelj's
formulas (that allow to find $\mu$ from the boundary values of $\hat
\mu$ on $(-1,1)$), that yield that the limit measure $\mu$ is
absolutely continuous and its density is
$$\mu'(x)=\frac{1}{\pi}\frac{dx}{\sqrt{1-x^2}}>0\,, \quad x\in(-1,1)\,,$$
(fact that we do not claim that was discovered by us; probably, it
is at least 150 years old).

In some more difficult situations the described elementary method
becomes technically involved, and we have to turn to other
resources. One of the advantages of the electrostatic interpretation
of the zeros of a sequence of polynomials is that it allows us to
guess what happens in the thermodynamic limit. For instance, if for
each $n$ the configuration of zeros is a global minimizer of the
total energy, then it is natural to expect (and can be proved) that
the asymptotics in a ``global sense'' should be well described (at
least, in the first approximation) by a continuous distribution
which is solution of the corresponding extremal problem for the
logarithmic potential energy in the class of all probability
measures (including the absolutely continuous ones). This
distribution is known as the equilibrium measure, probably, in an
external field (see \cite{Saff:97} for details and definitions).

Let us consider for example the generalized Lam\'{e} equation
(\ref{lame1}) under the assumptions
(\ref{simplezeros})--(\ref{lame2}). We are interested in the limit
$n\to \infty$ with the extra assumptions
\begin{equation*}%\label{rays}
  \lim_{n \to \infty} \frac{n_i}{n}=\theta_i\,, \quad i=1, \dots,
  p
\end{equation*}
(see Section \ref{sec3} for notation). Since for every $\vt{n}$ the
zeros minimize the total discrete energy of the system, we can
expect that any weak-* limit $\mu $ of the zero-counting unit
measures of the Heine-Stieltjes polynomials $E_{\vt{n}}$ solves a
similar, but continuous, extremal problem. Indeed, it was proved in
\cite{MS:01} that if $ \mathcal{M}=\mathcal{M}(\theta_1,
\dots,\theta_p)$ denotes the class of all probability measures $
\nu$ on $[-1, 1]$ such that
\[ \int_{a_{i-1}}^{a_i} d\nu= \theta_i, \quad i=1, \dots, p\,,
\]
then $\mu$ is the  minimizer of the logarithmic energy in the class
$ \mathcal{M}$. Existence and uniqueness of this minimizer is proved
by standard methods of potential theory. At this stage we can apply
the reduction to the Riccati form method described above to get the
full description of $\mu$. Let us summarize briefly some results
from \cite{MS:01}.

For any system of $p-1$ points
\begin{equation}\label{puntosBeta}
-1 \leq \beta_1 \leq \dots \leq \beta_{p-1} \leq 1
\end{equation}
we define the functions
$$
R(x):=\prod_{j=1}^{p-1} (x-\beta_j)\,,  \quad
H(x):=\sqrt{\frac{R(x)}{A(x)}} \sim \frac{1}{x}\,, \quad x \to
\infty\,.
$$
Then there exist $p-1$ points (\ref{puntosBeta}) uniquely determined
by the following system of equations:
$$
\Im \int_{a_{j-1}}^{a_j}
  H(x) \,dx=-\pi\,  \theta_j\,, \quad
  j=1, \dots, p-1\,,
$$
where we take the limit values of $H$ from the upper half plane. If
we introduce the counting function
\[
Z(x):=\left[\nu_A-\nu_R \right]\bigg( (-\infty,x]\bigg)\,,
\]
then $ \supp (\mu)=\overline{\{x \in \R:\, Z(x) =1 \}}$ and it
consists of at most $p-1$ disjoint intervals in $[-1,1]$.

Furthermore, $\mu$ is an absolutely continuous measure,
$$
\mu'(x)=-\frac{1}{\pi i}\, H(x)=\frac{1}{\pi }\, \left| H(x)
\right|\,, \quad x \in \supp (\mu)\,,
$$
and, for $z \notin \supp (\mu)$, $\hat \mu(z)=-H(z)$.

The methods described in this section yield many nice results, but
can fail for two reasons:
\begin{enumerate}
\item[(a)] If we have no a priori information on the location of the
zeros, then the reduction of the ODE to the Riccati form will give
us at most an expression of the Cauchy transform $\hat \mu$ of the
limit distribution $\mu$ in the domains (unknown) disjoint with
$\supp(\mu)$ (and $\supp(\mu)$ could be, eventually, a subset of
$\C$ of positive plane measure). Does this information determine
$\mu$? It is not clear, although for some specific expressions of
$\hat \mu$ this should be really the case.

\item[(b)] If we are analyzing Heine-Stieltjes polynomials whose
zeros provide a critical configuration for the total energy $E(X)$
(but not the global minimum!), then what measure $\mu$ should we
expect in the semiclassical limit? Clearly, $\mu$ is not necessarily
an equilibrium measure, but rather a \emph{critical measure}, that
on the real line can be characterized by the fact that its potential
(plus the external field, if exists) is constant on each connected
component of $\supp(\mu)$, but unlike in the equilibrium case, these
constants need not to be the same. However, in this description two
open problems remain: \emph{(i)} to prove that discrete critical
measures converge in the thermodynamic limit to continuous critical
measures, and that all continuous critical measures can be obtained
this way; \emph{(ii)} to find a feasible description of the
multi-parametric family of continuous critical measures in the given
class.
\end{enumerate}

\subsection*{Acknowledgements}

This research was supported, in part, by a grant of Direcci\'{o}n
General de Investigaci\'{o}n (Ministerio de Ciencia y Tecnolog\'{\i}a) of
Spain, project code BFM2003--06335--C03--02 (FM), a grant from the
Ministry of Education and Science of Spain, project code
MTM2005-08648-C02-01 (AMF, PMG), by Junta de Andaluc\'{\i}a, Grupo de
Investigaci\'{o}n FQM 0229 (AMF, PMG), by ``Research Network on
Constructive Complex Approximation (NeCCA)'', INTAS 03-51-6637 (FM,
AMF), and by NATO Collaborative Linkage Grant ``Orthogonal
Polynomials: Theory, Applications and Generalizations'', ref.
PST.CLG.979738 (AMF).

%%%%%%%%%%%%%%%%%%%%%%%%%%%%%%%%%%%%%%%%%%%%%%%%%%%%%%%%%%%%%%%%%%%%

\end{document}